\documentclass[11pt]{article}
\usepackage{amsmath}
\usepackage{mathrsfs}
\usepackage{amsfonts}

\usepackage{amsfonts, amsmath, amssymb}
\usepackage{amssymb,amsfonts,amsmath,
latexsym, epsfig,cite, psfrag,eepic,colordvi,color}
\usepackage{amscd,graphics}
\usepackage{lineno}
\allowdisplaybreaks

\newcommand{\qed}{\hfill $\Box $}
\newcommand{\pf}{\noindent {\bf Proof.} }

\newtheorem{theorem}{Theorem}[section]
\newtheorem{lemma}[theorem]{Lemma}

\newtheorem{conjecture}[theorem]{Conjecture}

\allowdisplaybreaks[4]
\usepackage{caption}
\begin{document}

\title{A proof of Frankl-Kupavskii's conjecture on edge-union condition \thanks{Supported by the National Natural
Science Foundation of China under grant No.12271425}}

\author{{Hongliang Lu and Xuechun Zhang}\\
\small School of Mathematics and Statistics\\
Xi'an Jiaotong University \\ \small  Xi'an, Shaanxi 710049, P.R.China \\ }

\date{}

\maketitle

\date{}

\maketitle {\small {\bfseries \centerline{Abstract}}

\vspace{3ex}

A  3-graph $\mathcal{F}$ is \emph{$U(s, 2s+1)$} if for any $s$ edges $e_1,...,e_s\in E(\mathcal{F})$, $|e_1\cup...\cup e_s|\leq 2s+1$. Frankl and  Kupavskii (2020) proposed the following conjecture: For any $3$-graph $\mathcal{F}$ with $n$ vertices, if $\mathcal{F}$ is $U(s, 2s+1)$, then
\[
e(\mathcal{F})\leq \max\left\{{n-1\choose 2}, (n-s-1){s+1\choose 2}+{s+1\choose 3}, {2s+1\choose 3}\right\}.
\]
 In this paper, we confirm Frankl and  Kupavskii's conjecture.

\vspace{3ex}

\section{Introduction}

A \emph{hypergraph} $\mathcal{F}$ consists of a vertex set $V (\mathcal{F})$ and an edge set $E(\mathcal{F})$, where every edge is a non-empty subset  of $V(\mathcal{F})$. A hypergraph $\mathcal{F}$ is
\emph{$k$-uniform}  if all its edges have size $k$ and we call it a $k$-graph for short. In particular, 2-graph  is called graph.
A \emph{matching} in a hypergraph $\mathcal{F}$ is a set of pairwise disjoint edges in $\mathcal{F}$.
 A  \emph{perfect matching} in a $k$-graph is a matching that covers every vertex. A \emph{maximum matching} in a $k$-graph is a matching of maximum size among all matchings in the graph.
 For a $k$-graph $\mathcal{F}$, we use $\nu(\mathcal{F})$ to
denote the  size of a maximum matching in $\mathcal{F}$.  Let
$G_1 = (V_1, E_1)$ and $G_2 = (V_2, E_2)$ be two graphs such that $V_1\cap V_2=\emptyset$. The \emph{union} of graphs $G_1$ and $G_2$ denoted by $G_1\cup G_2$   is the graph with vertex set $V_1\cup V_2$ and edge set $E_1\cup E_2$. The \emph{join} of graphs $G_1$ and $G_2$ denoted by $G_1\vee G_2$ is the graph union $G_1\cup G_2$ together with all the edges joining $V_1$ and $V_2$. For a graph $G$, we denote
the complement of $G$ by $\overline{G}$.


 We often identify $E(\mathcal{F})$ with $\mathcal{F}$ when there is no confusion and in particular, $|\mathcal{F}|$ denotes the number of
edges in $\mathcal{F}$.
 Given $T\subseteq V(\mathcal{F})$, let $\mathcal{F}-T$ denote the subgraph of
$\mathcal{F}$ with vertex set $V(\mathcal{F})\setminus T$ and edge set $\{e\in
\mathcal{F}\ : \ e\subseteq V(\mathcal{F})\setminus T\}$.
 Given a vertex set $X$ and integer $k\geq1$, let ${X\choose k}=\{S\subseteq X\ |\ |S|=k\}$. Let $n\geq 1$ be an integer and let $[n]=\{1,\ldots,n\}$. For a $k$-graph $\mathcal{F}$  and  $v\in\mathcal{F}$, let $N_{\mathcal{F}}(v)=\{e\in {V({\mathcal{F}})\choose k-1}\ |\ e\cup \{v\}\in E(\mathcal{F})\}$. When there is no confusion, we use  $N_{\mathcal{F}}(v)$ denote the $(k-1)$-graph with vertex set $V({\mathcal{F}})-v$ and edge set $N_{\mathcal{F}}(v)$.
  Let $m,i$ be two integers such that $1\leq i\leq  m\leq (n-2)/2$. Let $K^0_{i}$, $K^1_{2m-2i+1}$ and $K^2_{n-2m+i-1}$ denote complete graphs with vertex set $[i]$, $[2m-i+1]\setminus [i]$ and $[n]\setminus [2m-i+1]$ respectively and let $A^i_{n,m}:=K^0_{i}\vee (K^1_{2m-2i+1}\cup\overline{K^2_{n-2m+i-1}})$. Let $A^0_{n,m}:=K^1_{2m+1}\cup\overline{K^2_{n-2m-1}}$.

Let $p$, $r$ be  two positive integers such that $p\geq r$. Define
\begin{align*}
	\mathcal{A}_{p,r}:=\mathcal{A}(p,r,n,k):=\{A\in {[n] \choose k}:|A\cap [p]|\geq r \}.
\end{align*}
Let $k,s, q$ be three integers such that $2\leq s$ and $2\leq k\leq q<sk$. A $k$-graph $\mathcal{F}$ is said to have \emph{property   $U(s,q)$} if
\begin{align}
|e_1\cup...\cup e_s|\leq q
\end{align}
for any $s$ edges $e_1,...,e_s\in E(\mathcal{F})$. For shorthand, we will also say $\mathcal{F}$  is $U(s,q)$ to refer to this property.

Let
\begin{align*}
\mathcal{F}_1:=\mathcal{A}(1, 1, n, 3),\   \  \mathcal{F}_2:=\mathcal{A}(s+1,2,n,3),\  \
\mathcal{F}_3:=\mathcal{A}(2s+1,3,n,3),
\end{align*}
we can calculate $|\mathcal{F}_i|$, $i\in \{1,2,3\}$ from these definition,
\begin{align}\label{F1}
	|\mathcal{F}_1|={{n-1}\choose 2}=\frac{(n-1)(n-2)}{2},
\end{align}

\begin{align}\label{F2}	
	|\mathcal{F}_2|={{s+1}\choose 2}(n-s-1)+{{s+1}\choose 3}
	=-\frac{1}{3}s^3+(\frac{n}{2}-1)s^2+(\frac{n}{2}-\frac{2}{3})s,
\end{align}

\begin{align}\label{F3}
	|\mathcal{F}_3|={{2s+1}\choose 3}=\frac{4}{3}s^3-\frac{1}{3}s.
\end{align}
Frankl and Kupavskii \cite{FK21} conjectured that the maximum size of $|\mathcal{F}|$ when $k=3$ and $q=2s+1$.
\begin{conjecture}[Frankl and Kupavskii, \cite{FK21}]\label{conj}
Let $\mathcal{F}$ be a 3-graph with vertex set $[n]$. If $\mathcal{F}$ is $U(s,2s+1)$, then $|\mathcal{F}| \leq \max\{|\mathcal{F}_i| \ |\ i\in[3]\}$.
\end{conjecture}
In  this paper, we prove this conjecture.
\begin{theorem}\label{lu1}
Let $n,s$ be two integers such that $s\geq 2$ and $n\geq 2s+2$. Let $\mathcal{F}$ be a 3-graph with vertex set $[n]$. If $\mathcal{F}$ is $U(s,2s+1)$, then $|\mathcal{F}| \leq \max\{|\mathcal{F}_i| \ |\ i\in[3]\}$.

\end{theorem}

Let us recall some fundamental results of extremal set theory. Erd\H{o}s, Ko and Rado\cite{EKR} proved that if $\mathcal{F}\subseteq {n\choose k}$ satisfying $|F\cap F'|\geq t$ for any $F$, $F'\in \mathcal{F}$, then $|\mathcal{F}|\leq {n-t \choose {k-t}}$ for $n>n_0(k, t)$.
Moreover, they also  proved $n_0(k, 1)=2k$. We  call the conclusion as EKR Theorem in the following for short.
The smallest $n_0(k, t)$ of the EKR Theorem has been determined by Frankl \cite{EKRFK}} for
$t\geq 15$ and subsequently by Wilson \cite{EKRW} for all $t:n_0(k, t) = (k-t+1)(t+1)$. Some new generalizations and analogues of the EKR Theorem can be found \cite{DF,Fr91}. Frankl \cite{UC} added a union condition $U(s, n-t)$ to $\mathcal{F}$ and gave a upper bound to $|\mathcal{F}|$ when $t\leq {\frac{2^kk}{150}}$. Erd\H{o}s \cite{EMC} in 1965 made the following conjecture: If $\mathcal{F}\subseteq {n\choose k}$,  $\nu(\mathcal{F})\leq s$, then $|\mathcal{F}| \leq \max\{|\mathcal{A}_{s, 1}|, |\mathcal{A}_{(s+1)k-1, k}| \}$. It was one of the favourite problems of Erd\H{o}s and  attracted a lot of attention. 

The rest of the paper is organized as follows. In Section $2$, we first introduce some fundamental definitions and prove several lemmas needed in our proof.  
In Section $3$, we give the proof of Theorem \ref{lu1}. 

\section{Several Technical Lemma}

Let $k\geq 2$ be an integer.
Consider two sets $F_1=(a_1, a_2,\ldots, a_k)$ with $a_1<\cdots<a_k$ and $F_2=(b_1, b_2,\ldots, b_k)$ with $b_1<\cdots<b_k$. Then $F_1 \prec_s F_2$ iff $a_i\leq b_i$ for every $i\in [k]$. We say that a family $\mathcal{F}\subseteq{[n]\choose k}$ is \emph{stable} if $F_2\in \mathcal{F}$ and $F_1\prec_s F_2$ implies $F_1\in \mathcal{F}$. 

\begin{lemma}\label{S-subgraph}
Let $G$ be a stable graph with vertex set $[n]$ and $\nu(G)=m$. Let $r:=\max \{i\ |\ \nu(G-\{0,\ldots,i\})=\nu(G)-i\}$. Then $G$ is a subgraph of $A_{n,m}^r$. 
\end{lemma}

\pf Let $\mathcal{M}=\{\{j,2m+3-j\}\ |\ 1\leq j\leq m+1\}$. Since $\nu(G)=m$, then $\mathcal{M}$ is not a matching of $G$. So there exists $j\in [m+1]$ such that $\{j,2m+3-j\}\notin G$. Since $G$ is stable, for any edge $\{u,v\}\in G$ with $u<v$, either $u<j$ or $v<2m+3-j$. Thus
$G$ is a subgraph of $A_{n,m}^{j-1}$. We choose maximum $j$ such that $G$ is a subgraph of $A_{n,m}^{j-1}$.

Next we show that $r=j-1$. Note that
\[
m+1-j=\nu(G)+1-j\leq \nu(G-[j-1])\leq \nu(A_{n,m}^{j-1}-[j-1])= m+1-j.
\]
 Thus we have $\nu(G-[j-1])=\nu(G)-j+1$. Thus by the definition of $r$, we have $r\geq j-1$.
 We claim  $\{\{i+1,2m+2-i\} |\ j\leq i\leq m\}$ is a matching of $G$. Otherwise, there exists $t\in [m]-[j-1]$ such that $\{t+1,2m+2-t\}\notin G$.
  Hence for any  $\{a,b\}\in G$ with $a<b$, either $a\leq t$ or $b\leq 2m+1-t$. So $G$ is a subgraph of $A^{t}_{n,m}$, contradicting  the choice of $j$ since $j-1<t$. This completes the proof. \qed
 %
%
%

The following result can be found in \cite{UC,AF85}.
\begin{lemma}\label{S-matching}
Let $G$ be a stable graph with vertex set $[n]$ and $\nu(G)=s$. Then $\{\{i,2s-i+1\}\ |\ 1\leq i\leq s\}$ is a matching of size $s$ in $G$.
\end{lemma}


A basic result of Frankl \cite{ST,shifted} allows us to restrict our proof  to stable family.
\begin{lemma}\label{stable}
If the family $\mathcal{F}\subseteq {[n]\choose k} $
is $U(s, q)$, then there exists a stable and $U(s, q)$ family $\mathcal{F}'\subseteq {[n]\choose k} $
such that $|\mathcal{F}|= |\mathcal{F}'|$.
\end{lemma}
By lemma \ref{stable},   we may assume that $\mathcal{F}$ is stable in our proof.
We denote the shadow of $\mathcal{F}$  by $\partial \mathcal{F}$.
\begin{align}
	\partial \mathcal{F}=\{G\in {[n]\choose {k-1}}:\exists  F\in \mathcal{F}, G\subset F\}
\end{align}
When $k=3$, $\partial \mathcal{F}\subseteq {V(\mathcal{F})\choose 2}$ and so $\partial \mathcal{F}$ can be seen as a graph. Note that if $\mathcal{F}$ is stable, then $\partial \mathcal{F}$ is stable.

\section{Proof of Theorem \ref{lu1}}

If $\nu(\mathcal{F})=1$, then by EKR Theorem,  $|\mathcal{F}|\leq |\mathcal{F}_1|\leq \max\{|\mathcal{F}_i|\ |\ i\in [3]\}$. So we may assume that $\nu(\mathcal{F})\geq 2$.
Note that    $\nu(\mathcal{F})\geq 2$ implies that $\mathcal{F}\subseteq {[n] \choose 3}$ is not $U(2, 5)$. So we have the following result.
\begin{lemma}\label{s=2}
Let $n$ be an integer such that $n\geq 7$, If $\mathcal{F}\subseteq {[n] \choose 3}$ is $U(2, 5)$, then  $|\mathcal{F}| \leq \max \{|\mathcal{F}_i|\ |\ i\in[3]\}$.
\end{lemma}

Thus, we can restrict our proof to $s\geq 3$. By Lemma \ref{stable},  we give the proof of  Theorem \ref{lu1} by proving Lemmas \ref{m_lem1}, \ref{m_lem2} and \ref{m_lem3}.

\begin{lemma}\label{m_lem1}
	Let $n,s$ be two integer such that $s\geq 3$ and $n\geq 2s+4$. Let $\mathcal{F}$ be a stable 3-graph with vertex set $[n]$. 	If $\mathcal{F}$ is $U(s,2s+1)$, $\nu(\mathcal{F})\geq 2$ and $\nu(\partial\mathcal{F})\geq s+2$, then   $|\mathcal{F}| \leq\max \{|\mathcal{F}_i|\ |\ i\in[3]\}$.
\end{lemma}
\noindent\textbf{Proof.}  Since $\mathcal{F}$ is stable and $\nu(\mathcal{F})\geq 2$, we have $\{2,3,4\}\in \mathcal{F}$.    We first prove the following claim.

 \medskip
\textbf{Claim 1.~} $\nu(\partial\mathcal{F}-[4])\leq s-2$.

Otherwise, suppose that $\nu(\partial\mathcal{F}-[4])\geq s-1$. Then $\partial\mathcal{F}-[4]$ has a matching of size $s-1$. Let $M$ be a matching of size $s-1$ in $\partial\mathcal{F}-[4]$. Write $M':=\{e\cup \{1\}\ |\ e\in M\}\cup\{\{2,3,4\}\}$. Since $\mathcal{F}$ is stable, we have  $M'\subseteq \mathcal{F}$. One can see that $|M'|=s$ and  $|\cup_{e\in M'}e|=2s+2$, a contradiction since $\mathcal{F}$ is $U(s,2s+1)$. This completes the proof of Claim 1.

  Since $\nu(\partial\mathcal{F})\geq s+2$, by Claim 1, we have $\nu(\partial\mathcal{F})=s+2$ and $\nu(\partial\mathcal{F}-[4])=s-2$.

 \medskip
\textbf{Claim 2.~} If $M_0$ is a matching of size $s$ in $\partial\mathcal{F} -\{1,2\}$, then $ N_{\mathcal{F}}(2)\cap M_0=\emptyset$.

 Otherwise, suppose that $ N_{\mathcal{F}}(2)\cap M_0\neq \emptyset$. Let $e\in M_0\cap N_{\mathcal{F}}(2)$ and $M_0'=\{f\cup \{1\}\ |\ f\in M_0\setminus e\}\cup \{e\cup \{2\}\}$. By the stability of $\mathcal{F}$, one can see that $M_0'\subseteq \mathcal{F}$ and $|V(M_0')|=2s+2$. So $\mathcal{F}$ is not $U(s,2s+1)$, a contradiction.  This completes the proof of Claim 2.

 \medskip
\textbf{Claim 3.~} $e=\{2,3,2s+2\}\notin \mathcal{F}$.

By contradiction. Suppose that $\{2,3,2s+2\}\in \mathcal{F}$. By Claim 1, $\partial\mathcal{F}-[3]$ has a matching of size $s-1$. By Lemma \ref{S-matching}, $M_1:=\{\{i,2s+5-i\}\ |\ 4\leq i\leq s+2\}$ is matching of size $s-1$ in $\partial\mathcal{F}-[3]$.
Let $M_1'=\{f\cup\{1\}\ |\ f\in M_1\}\cup \{\{2,3,2s+2\}\}$. Then we have $|M_1'|=s$, $|V(M_1')|=2s+2$. Since $M_1\subseteq  N_{\mathcal{F}}(1)$, we have $M_1'\subseteq \mathcal{F}$. So $\mathcal{F}$ is not $U(s,2s+1)$, a contradiction.  This completes the proof of Claim 3.

 Let $r=\max\{i\ |\ \nu(\partial\mathcal{F}-[i])=\nu(\partial\mathcal{F})-i\}$. From Claim 1 we have $r\geq 4$.
By Lemma \ref{S-subgraph}, $\partial\mathcal{F}$ is a subgraph of $K_{r}^0\vee (K_{2s-2r+5}^1\cup K_{n-2s+r-5}^2)$.
By Lemma \ref{S-matching}, $\{\{i,2s+5-i\}\ |\ 3\leq i\leq s+2\}$ is a matching of size $s$ in $\partial\mathcal{F}- \{1, 2\}$.
Note that $\{s+2,s+3\}\in \partial\mathcal{F} $ and  by Claim 3, we have $\{s+2,s+3,2\}\notin \mathcal{F}$. By the stability of $\mathcal{F}$, every edge  of $\mathcal{F} $ containing no $1$ intersects $[s+1]$ at least two vertices, i.e.,  $\mathcal{F}-\{1\}$ is a subgraph of $\mathcal{F}_2-\{1\}$.
Thus we have
\begin{align*}
|\mathcal{F}\setminus \mathcal{F}_2|&=|N_{\mathcal{F}}(1)\setminus N_{\mathcal{F}_2}(1)|\\
&\leq \left|\left(\left({[n]-[1]\choose 2}\setminus{[n]-[r]\choose 2}\right)\cup {[2s-r+5]-[r]\choose 2}\right)\setminus\left({[n]-[1]\choose 2}\setminus{[n]-[1+s]\choose 2}\right)\right|\\
&={s-r+4\choose 2}\\
&\leq {s\choose 2}\quad \mbox{(since $r\geq 4$)}.
\end{align*}

Note that $\mathcal{M}:=\{\{i+2,2s+3-i\}\ |\ 1\leq i\leq s\}$ is a matching of size $s$ in $\partial\mathcal{F}-[2]$. By Claim 2, we have $\mathcal{M}\cap N_{\mathcal{F}}(2)=\emptyset$.

 \medskip
\textbf{Claim 4.~} For any two distinct edges $e_1,e_2\in \mathcal{M}$, ${e_1\cup e_2\choose 3}\cap \mathcal{F}= \emptyset$.

Otherwise, suppose that ${e_1\cup e_2\choose 3}\cap \mathcal{F}\neq \emptyset$. Let $e\in {e_1\cup e_2\choose 3}\cap \mathcal{F}$ and let $v\in (e_1\cup e_2)\backslash e$. Note that $v\leq 2s+2$. So we have $\{1,2,v\}\in \mathcal{F}$.
Let $M_2=\{e,\{1,2,v\}\}\cup \{f\cup \{1\}\ |\ f\in \mathcal{M}\backslash \{e_1,e_2\}\}$. Then we have $|M_2|=s$ and $|V(M_2)|=2s+2$. By the stability, One can see that $M_2\subseteq \mathcal{F}$, contradicting that $\mathcal{F}$ is $U(s,2s+1)$. This completes the proof of Claim 4.

Write $A=\{\{s+2,s+3\}, \{3,2s+2\}\}$.  Recall that $\mathcal{F}$ and $\partial\mathcal{F}$ are stable. By Claims 2 and 4, $\{i\}\cup e\notin \mathcal{F}$ for any $i\in \{2,3\}$ and $e\in \mathcal{M}\backslash A$. By Claim 3, $\{j\}\cup e\notin \mathcal{F}$ for any $j\in [n]-[2s+1]$ and $e\in {[s+1]-\{1\}\choose 2}$.
Hence by Claim 4,   we have
\begin{align*}
|\mathcal{F}_2\setminus \mathcal{F}|&\geq\left|{\{2,3\}\choose 1}\right|\left|\mathcal{M}\backslash A\right|+{s\choose 2}(n-2s-1)+{4\choose 3}\left|{\mathcal{M}\backslash A\choose 2}\right|\\
&=2(s-2)+{s\choose 2}(n-2s-1)+4{s-2\choose 2}\\
&=2(s-2)^2+{s\choose 2}(n-2s-1)\\
&>{s\choose 2}\quad\mbox{(since $n\geq 2s+4$)}.
\end{align*}
So we have $|\mathcal{F}_2\setminus \mathcal{F}|>|\mathcal{F}\setminus \mathcal{F}_2|$, which implies that $|\mathcal{F}_2|>|\mathcal{F}|$. This complete the proof. \qed

\begin{lemma}\label{m_lem2}
	Let $n,s$ be two integer such that $s\geq 3$ and $n\geq 2s+2$. Let $\mathcal{F}$ be a stable 3-graph with vertex set $[n]$. 	If $\mathcal{F}$ is $U(s,2s+1)$ and $\nu(\partial\mathcal{F})\leq s$, then $|{\mathcal{F}}| \leq \max \{|\mathcal{F}_i|\ |\ i\in[3]\}$.
\end{lemma}

\noindent\textbf{Proof.}  Without loss of generality, we may assume that $\nu(\partial\mathcal{F})=s$ since we may add some edges containing $1$ to $\mathcal{F}$ such that $\nu(\partial\mathcal{F})=s$. Let $r:=\max \{i\ |\ \nu(\partial\mathcal{F}-\{0,\ldots,i\})=\nu(\partial\mathcal{F})-i\}$.

Since $\partial\mathcal{F}$ is stable and $\nu(\partial\mathcal{F})=s$, by Lemma \ref{S-subgraph}, $\partial\mathcal{F}$ is a subgraph of $A^r_{n,s}$.
So for every edge $e$ of $\mathcal{F}$, we have $|e\cap [r]|=2$ or $e\subseteq [2s-r+1]$. Thus  we have
\begin{align*}
	|\mathcal{F}|&\leq {r\choose 2}(n-2s+r-1)+{{2s-r+1}\choose 3}\\
	&=\frac{1}{3}r^3+(\frac{n}{2}-1)r^2-(\frac{n}{2}-s+2s^2-\frac{2}{3})r+\frac{4}{3}s^3-\frac{1}{3}s.
\end{align*}
Let $f(r)=\frac{1}{3}r^3+(\frac{n}{2}-1)r^2-(\frac{n}{2}-s+2s^2-\frac{2}{3})r+\frac{4}{3}s^3-\frac{1}{3}s$. Then $f(r)$ is a convex function on interval $[0,s]$ since $f''(r)=2r+n-2\geq 0$ for $0\leq r\leq s$. Thus
\begin{align*}
	|\mathcal{F}|&\leq \max\{f(0), f(s)\}\\
   &=\max\{{2s+1\choose 3}, {s\choose 2}(n-s)+{s\choose 3} \}\\
   &\leq \max\{|\mathcal{F}_2|,|\mathcal{F}_3|\}.
\end{align*}
This completes the proof. \qed


\begin{lemma}\label{<=4}
Let $\mathcal{F}\subseteq {[6]\choose 3}$ such that $e\cap \{i,7-i\}\neq \emptyset$ for all $e\in\mathcal{F}$, $1\leq i\leq 3$. If $\nu(\mathcal{F})=1$, then $|\mathcal{F}|\leq 4$.

\end{lemma}

\pf Let $M_1=\{\{1,2,3\},\{4,5,6\}\}$, $M_2=\{\{1,2,4\},\{3,5,6\}\}$, $M_3=\{\{1,4,5\},\{2,3,6\}\}$ and $M_4=\{\{1,3,4\},\{2,5,6\}\}$. By the definition of $\mathcal{F}$, we have $\mathcal{F}\subseteq \cup_{i=1}^4M_i$.
Since $\nu(\mathcal{F})=1$, then $|M_i\cap \mathcal{F}|\leq 1$ for $1\leq i\leq 4$. Thus $|\mathcal{F}|\leq \sum_{i=1}^4|M_i\cap \mathcal{F}|\leq 4$. This completes the proof. \qed

\begin{lemma}\label{m_lem3}
	Let $n,s$ be two integers such that $s\geq 3$ and $n\geq 2s+2$. Let $\mathcal{F}$ be a stable 3-graph with vertex set $[n]$. 	If $\mathcal{F}$ is $U(s,2s+1)$ and $\nu(\partial\mathcal{F})= s+1$, then $|{\mathcal{F}}| \leq \max\{|\mathcal{F}_i|\ |\ i\in[3]\}$.
\end{lemma}
\noindent\textbf{Proof.} Let $r:=\max \{i\ |\ \nu(\partial\mathcal{F}-\{0,\ldots,i\})=\nu(\partial\mathcal{F})-i\}$. 
By  Lemma \ref{S-subgraph}, $\partial\mathcal{F}$ is a subgraph of $A_{n,s}^r$.
Now we discuss  four cases.

%

 \textbf{Case 1.} $3\leq r\leq s+1$.

For every edge $e\in \mathcal{F}$, $e\subseteq [2s+3-r]$ or $|e\cap [r]|=2$. So we have
 \begin{align*}
 	|\mathcal{F}|&\leq {r\choose 2}(n-2s+r-3)+{{2s-r+3}\choose 3}\\
 	&=\frac{1}{3}r^3+(\frac{n}{2}-1)r^2-(\frac{n}{2}+2s^2+3s+\frac{1}{3})r+\frac{4}{3}s^3+4s^2+\frac{11}{3}s+1,
 \end{align*}
 Let $g(r):=\frac{1}{3}r^3+(\frac{n}{2}-1)r^2-(\frac{n}{2}+2s^2+3s+\frac{1}{3})r+\frac{4}{3}s^3+4s^2+\frac{11}{3}s+1.$ Note that
$g(r)$ is a convex function on interval $[3,s+1]$ since $g''(r)=2r+n-2\geq 0$ for $3\leq r\leq s+1$. Thus
\begin{align*}
	|\mathcal{F}|&\leq \max\{g(3), g(s+1)\}\\
	&=\max\{{2s\choose 3}+3(n-2s),  {s+2\choose 3}+{s+1\choose 2}(n-s-2)\}\\
&=\max\{{2s\choose 3}+3(n-2s),  {s+1\choose 3}+{s+1\choose 2}(n-s-1)\}\\
&=\max\{3n+\frac{4}{3}s^3-2s^2-\frac{16}{3}s, |\mathcal{F}_2|\}.
\end{align*}
Next we discuss two subcases.

\medskip
\textbf{Subcase 1.1.~$n\leq (2s^2+5s)/3$.}

 From (\ref{F3}) we have
\begin{align*}
   	|\mathcal{F}_3|-g(3)&=-3n+2s^2+5s\geq 0.   	
\end{align*}
Hence, $g(3)\leq |\mathcal{F}_3|$. On the other hand, $g(s+1)=|\mathcal{F}_2|$. So we have	$|\mathcal{F}|\leq \max\{g(3), g(s+1)
\}\leq \max\{|\mathcal{F}_i|\ |\ i\in[3]\}$.

\medskip
\textbf{Subcase 1.2.~$n> (2s^2+5s)/3$.}

Note that $n> (2s^2+5s)/3>10s/3$. Thus we have
\begin{align*}
g(s+1)-g(3)&=-\frac{5}{3}s^3+(\frac{n}{2}+1)s^2+(\frac{n}{2}+\frac{14}{3})s-3n\\
&>-\frac{5}{3}s^3+(\frac{5}{3}s+1)s^2+(\frac{5}{3}s+\frac{14}{3})s-10s\quad \mbox{(since $s\geq 3$)}\\
&=\frac{8s}{3}(s-2)>0.
\end{align*}
Thus by (\ref{F2}), one can see that
\begin{align*}
	|\mathcal{F}|&\leq \max \{g(3), g(s+1)\}=g(s+1)=|\mathcal{F}_2|.
\end{align*}

\medskip
 \textbf{Case 2.} $r=0$.

Then $\partial\mathcal{F}$ is a subgraph of $K_{2s+3}$. Thus $\mathcal{F}\subseteq {[2s+3]\choose 3}$. From the definition of $r$  we have $\nu(\partial\mathcal{F}-\{1\})=\nu(\partial\mathcal{F})=s+1$, by Lemma \ref{S-matching}, $\mathcal{M}:=\{\{i+1,2s+4-i\ |\ 1\leq i\leq s+1\}\}$ is a matching of size $s+1$ in $\partial\mathcal{F}$.

\medskip
\textbf{Claim 1.} For any $e_1,e_2\in \mathcal{M}$, ${e_1\cup e_2\choose 3}\cap \mathcal{F}=\emptyset$.

Otherwise, suppose that ${e_1\cup e_2\choose 3}\cap \mathcal{F}\neq\emptyset$. Let $f\in {e_1\cup e_2\choose 3}\cap \mathcal{F}$. Write $M=\{f\}\cup \{e\cup \{1\}\ |\ e\in \mathcal{M}-\{e_1,e_2\}\}$. One can see that $M\subseteq \mathcal{F}$ and $|V(M)|=2s+2$, contradicting that $\mathcal{F}$ is $U(s,2s+1)$. This complete the proof of Claim 1.

\textbf{Claim 2.} For any three distinct edges $e_1,e_2,e_3\in \mathcal{M}$, $\mathcal{F}[\cup_{i\in [3]}e_i\cup\{1\}]$ contains no  matchings of size two. 

Otherwise, there exists $e_1,e_2,e_3\in \mathcal{M}$ such that  $\mathcal{F}[(\cup_{i\in [3]}e_i)\cup \{1\} ]$ contains a  matching of size two, say $\{f_1,f_2\}$.  Write $M:=\{f_1,f_2\}\cup \{e\cup \{1\}\ |\ e\in \mathcal{M}-\{e_i\ |\ i\in [3]\}\}$. One can see that  $M\subseteq \mathcal{F}$ and $|V(M)|\geq 2s+2$, contradicting that $\mathcal{F}$ is $U(s,2s+1)$. This complete the proof of Claim 2.


By Lemma \ref{<=4}, Claims 1 and 2, we have
\begin{align*}
|\mathcal{F}|-|\mathcal{F}_3|&\leq {2s+3\choose 3}-\left|{\mathcal{M}\choose 2}\right|{4\choose 3}-4\left|{\mathcal{M}\choose 3}\right|-{2s+1\choose 3}\\
&={2s+2\choose 2}+{2s+1\choose 2}-4{s+1\choose 2}-4{s+1\choose 3}\\
&=-\frac{2}{3}s^3+2s^2+\frac{8}{3}s\\
&\leq -\frac{2}{3}s^2+\frac{8}{3}s\leq 0\quad \mbox{(when $s\geq 4$)},
\end{align*}
i.e., $|\mathcal{F}|\leq |\mathcal{F}_3|$.
Next we consider $s=3$. Then $\mathcal{M}=\{\{i+1,10-i\}\ |\ i\in [4]\}$. If $\{1,6,7\}\in \mathcal{F}$, by Claim 2, $\{2,4,5\}\notin \mathcal{F}$. By the stability of $\mathcal{F}$, we have $\{2,i,j\}\notin \mathcal{F}$, $i,j\in \{4,5,\ldots,9\}$, thus
\begin{align*}
|\mathcal{F}|&\leq {9\choose 3}-\left|{[9]-[2]\choose 3}\right|-\left|{[9]-[3]\choose 2}\right|=34<35=|\mathcal{F}_3|.
\end{align*}

So we  may  assume that $\{1,6,7\}\notin \mathcal{F}$. By Claim 1, $\{2,5,6\},\{2,4,7\}\notin \mathcal{F}$, and so by stability, we have $\{2,5,7\},\{2,5,8\},\{2,4,8\}\notin \mathcal{F}$. Hence
\begin{align*}
|\mathcal{F}|&\leq {9\choose 3}-\left|{[5]\choose 1}\right| \left|{[9]-[5]\choose 2}\right|-2\left|{\mathcal{M}\choose 2}\right|-\left|{[9]-[5]\choose 3}\right|-3\\
&=35=|\mathcal{F}_3|.
\end{align*}

 \textbf{Case 3.} $r=1$.

By Lemma \ref{S-subgraph}, $\partial \mathcal{F}$ is a subgraph of $A_{n,s+1}^1$. So we have $\mathcal{F}\subseteq {[2s+2]\choose 3}$.
By Lemma \ref{S-matching}, $\mathcal{M}:=\{\{i,2s+3-i\}\ |\ 1\leq i \leq s+1\}$ is a matching of size $s+1$ of $\partial \mathcal{F}$ . If $s\geq 4$, by Claim 2 and Lemma \ref{<=4}, we have
\begin{align*}
|\mathcal{F}|-|\mathcal{F}_3|&\leq {2s+2\choose 3}-4\left|{\mathcal{M}\choose 3}\right|-{2s+1\choose 3}\\
&=-\frac{2}{3}s^3+2s^2+\frac{5}{3}s\\
&\leq -\frac{1}{3}(2s^2-5s)\quad\mbox{(since $s\geq 4$)}\\
&<0.
\end{align*}
So we may assume that $s=3$.
Then  $\mathcal{M}:=\{\{i,9-i\}\ |\ 1\leq i \leq 4\}$. If $\{1,2,8\}\notin \mathcal{F}$, then we have $\mathcal{F}\subseteq {[2s+1]\choose 3}$ and so $|\mathcal{F}|\leq |\mathcal{F}_3|$. Hence we may assume $\{1,2,8\}\in \mathcal{F}$. By Claim 2 and the stability of  $\mathcal{F}$, we have $\{4,5,7\},\{4,5,8\},\{3,6,7\},\{3,6,8\}\notin \mathcal{F}$. Moreover, one can see that $|\{\{1,7,8\},\{2,3,6\}\}\cap\mathcal{F}|=1$ by Claim 2. Let $e\in \{\{1,7,8\},\{2,3,6\}\}$ such that $e\notin \mathcal{F}$. Write
\[
M=\{e,\{4,5,7\},\{4,5,8\},\{3,6,7\},\{3,6,8\}\}.
\]
Thus we have
\begin{align*}
|\mathcal{F}|-|\mathcal{F}_3|&\leq {8\choose 3}-{7\choose 3}-4{|\mathcal{M}|\choose 3}-|M|=0,
\end{align*}
and so $|\mathcal{F}|\leq |\mathcal{F}_3|$.

\medskip
 \textbf{Case 4.} $r=2$.

By Lemma \ref{S-subgraph}, $\partial \mathcal{F}$ is a subgraph $A_{n,s+1}^2$. By the definition of $r$, we have $\partial\mathcal{F}-[3]$ has a matching size $s-1$. By Lemma \ref{S-matching},  $\mathcal{M}:=\{\{i,2s+5-i\}\ |\ 4\leq i\leq s+2\}$ is a matching of $\partial\mathcal{F}-[3]$. If $\{1,2,2s+2\}\notin \mathcal{F}$, then $\partial \mathcal{F}\subseteq {[2s+1]\choose 2}$ and $\nu(\mathcal{F})\leq s$, a contradiction. So we may assume that $\{1,2,2s+2\}\in \mathcal{F}$.


\medskip
\textbf{Claim 3.}  $\mathcal{M}\cap N_{\mathcal{F}}(3)=\emptyset$. In particular, $\{3,s+2,s+3\}\notin \mathcal{F}$.

Otherwise, suppose that $e\in\mathcal{M}\cap N_{\mathcal{F}}(3)$.  Let $M:=\{\{3\}\cup e,\{1,2,2s+2\}\}\cup\{f\cup\{1\}\ |\  f\in \mathcal{M}\setminus e\}$. Since $\mathcal{F}$ is stable, $M\subseteq \mathcal{F}$. Moreover, $|M|=s$ and $|V(M)|=2s+2$, contradicting that $\mathcal{F}$ is $U(s,2s+1)$.


Write $A=\{3,\ldots,s+1\}$ and $B=[2s]-[s+1]$. By Claim 3, we have $\{3,4,2s+1\}\notin \mathcal{F}$. So for any $\{a,b\}\in {A\choose 2}$, $\{a,b,2s+1\}\notin \mathcal{F}$.
Then  by Claim 3, we have
\begin{align*}
|\mathcal{F}|&\leq {2s+1\choose 3}+(n-2s-1)-{|A|\choose 1}{|B|\choose 2}-{|A\cup B|\choose 2}\\
&={2s+1\choose 3}+(n-2s-1)-(s-1){s-1\choose 2}-{2s-2\choose 2}\\
&={2s+1\choose 3}+n-3-\frac{1}{2}(s-1)(s^2+s),
\end{align*}
i.e.,
\begin{align}\label{n_s_bound}
|\mathcal{F}|\leq{2s+1\choose 3}+n-3-\frac{1}{2}(s-1)(s^2+s).
\end{align}
Consider $n\leq \frac{1}{2}(s-1)(s^2+s)+3$. Then we have $|\mathcal{F}|\leq |\mathcal{F}_3|$. So we may assume that $n> \frac{1}{2}(s-1)(s^2+s)+3\geq 3s+6$ as $s\geq 3$. Then by (\ref{n_s_bound}),
\begin{align*}
|\mathcal{F}|-|\mathcal{F}_2|&\leq {2s+1\choose 3}+n-\frac{1}{2}(s-1)(s^2+s)-4-{s+1\choose 2}(n-s-1)-{s+1\choose 3}\\
&\leq {2s+1\choose 3}+3s+7-\frac{1}{2}(s-1)(s^2+s)-4-{s+1\choose 2}(2s+6)-{s+1\choose 3}\quad\mbox{(since $n\geq 3s+7$)}\\
&=\frac{1}{3}s(4s^2-1)-\frac{1}{2}(s+1)(3s^2+5s-6)-\frac{1}{6}s(s^2-1)\\
&\leq \frac{1}{3}s(4s^2-1)-\frac{3}{2}s^3-\frac{1}{6}s(s^2-1)\quad\mbox{(since $s\geq 3$)}\\
&=-\frac{1}{3}s^3-\frac{1}{6}s<0,
\end{align*}
i.e., $|\mathcal{F}|\leq|\mathcal{F}_2|$. This completes the proof. \qed

%
%
%



{\bf Acknowledgment.} The authors would like to thank Professor Peter Frankl  for his  helpful suggestions and comments.

\end{document}